\documentclass[preprint,review,12pt]{elsarticle}
\usepackage{amsfonts}
\usepackage{amssymb}
\usepackage{graphicx}
\usepackage[top=1.2in, bottom=1.2in, left=1.2in, right=1.0in]{geometry}
\newtheorem{Theorem}{Theorem}[section]

\newtheorem{Remark}{Remark}[section]

\newcommand{\be}{\begin{equation}}
\newcommand{\ee}{\end{equation}}
\newcommand\bes{\begin{eqnarray}}
\newcommand\ees{\end{eqnarray}}
\newcommand{\bess}{\begin{eqnarray*}}
\newcommand{\eess}{\end{eqnarray*}}

\usepackage{amssymb}

\journal{Physic A}

\begin{document}

\begin{frontmatter}



\title{Cross-diffusion induced Turing instability in two-prey one-predator system\tnoteref{label1}}
\tnotetext[label1]{The work supported by PRC grant NSFC (61472343, 11201406) and NSF of YZU (2014CXJ002).}
\author{Yahui Chen$^a$, Canrong Tian$^b$}

\author{Zhi Ling$^a$}
\ead{zhling@yzu.edu.cn}


\address{$^a$School of Mathematical Science, Yangzhou University, \\ Yangzhou 225002, P.R. China}
\address{$^b$Department of
Basic Science, Yancheng Institute of Technology,\\ Yancheng 224003,
P. R. China.}

\begin{abstract}
In this paper, we study a strongly coupled two-prey one-predator
system. We first prove the unique positive equilibrium solution is
globally asymptotically stable for the corresponding kinetic system
(the system without diffusion) and remains locally linearly stable
for the reaction-diffusion system without cross-diffusion, hence it
does not belong to the classical Turing instability scheme. Moreover
we prove that the positive equilibrium solution is globally
asymptotically stable for the reaction-diffusion system without
cross-diffusion. But it becomes linear unstable only when
cross-diffusion also plays a role in the reaction-diffusion system,
thus it is a cross-diffusion induced instability. Finally, the
corresponding numerical simulations are also demonstrated and we
obtain the spatial patterns.
\end{abstract}


\end{frontmatter}


\section{Introduction}
In 2009, Elettreby considered the following prey-predator model
\cite{Ele}
\begin{eqnarray}
\label{eq0}\left\{
\begin{array}{l}
\displaystyle u_{1}'=au_1(1-u_1)-u_1u_3:=u_1f_1(u_1,u_3),\\
\displaystyle u_{2}'=bu_2(1-u_2)-u_2u_3:=u_2f_2(u_2,u_3),\\
\displaystyle u_{3}'=-cu_3^2+(du_1+eu_2)u_3:=u_3f_3(u_1,u_2,u_3),
\end{array} \right.
\end{eqnarray}
where $u_{1}, u_{2}$ and $u_{3}$ are the population densities of
three species. This system models the dynamic of two-prey
one-predator ecosystem, i.e. the third species preys on the second
and the first one. In the absence of any predation, each term of
preys grows logistically. The effect of the predation is to reduce
the prey growth rate. In the absence of any prey for sustenance, the
predator's death rate results in inverse decay, which is the term
$-cu_{3}^{2}$. The prey's contribution to growth rate of the
predators are respectively $du_{1}u_{3}$ and $eu_{2}u_{3}$. They
studyed the global stability and persistence of the model.

However, in reality, individual organisms are distributed in space.
We can use the reaction-diffusion equations to establish
spatio-temporal dynamical system which can model the pursuit-evasion
phenomenon (predators pursuing prey and prey escaping predators) in
the prey-predator system. Therefore, in present paper we further
investigate the following reaction-diffusion model with
cross-diffusion:
\begin{eqnarray}
\label{eq1}\left\{
\begin{array}{lll}
u_{1t}-\Delta[(k_{11}+k_{13}u_3) u_1]=au_1(1-u_1)-u_1u_3, &  \mbox{in}\ \Omega\times(0,\infty),\\
u_{2t}-\Delta[(k_{22}+k_{23}u_3) u_2]=bu_2(1-u_2)-u_2u_3, &  \mbox{in}\  \Omega\times(0,\infty),\\
u_{3t}-\Delta[(k_{31}u_1+k_{32}u_2+k_{33})u_3]=-cu_3^2+(du_1+eu_2)u_3, & \mbox{in}\  \Omega\times(0,\infty),\\
\frac{\partial u_1}{\partial \eta}=\frac{\partial u_2}{\partial
\eta}=\frac{\partial u_3}{\partial \eta}=0, &  \mbox{on}
\ \partial \Omega \times (0,\infty),\\
u_1(x,0)=u_{10}(x), u_2(x,0)=u_{20}(x), u_3(x,0)=u_{30}(x), &
\mbox{in}\ \ \Omega, \end{array} \right.
\end{eqnarray}
where $\Omega$ is a bounded domain in $R^N$ with smooth boundary
$\partial\Omega$. $\eta$ is the unit outward normal to
$\partial\Omega$. The homogeneous Neumann boundary condition
indicates that there is zero population flux across the boundary.
The parameters $a, b, c, d, e$ and $k_{ij}$ $(1\leq i$, $j\leq3)$
are all positive constants. $k_{ii}$ is the diffusion rate of $i$-th
species. This diffusion term represent simple Brownian type motion
of particle dispersal. $k_{ij}$ $(i\neq j)$ is the cross-diffusion
rate of $i$-th species. It is necessary to note that the
cross-diffusion coefficient may be positive or negative. The
positive cross-diffusion coefficient represents that one species
tends to move in the direction of lower concentration of another
species. On the contrary, the negative cross-diffusion coefficient
denotes the population flux of one species in the direction of
higher concentration of another species. Here the cross-diffusion
term presents the tendency of predators to avoid the group defense
by a large number of prey, i.e. the predator diffuses in the
direction of lower concentration of the prey species. More
biological background can be found in \cite{ Can,lou,Okubo}.

As we know, the problem of cross-diffusion was proposed first by
Kerner \cite{K1} and first applied to competitive population systems
by Shigesada et al. \cite{SKT}. Since then, the role of
cross-diffusion in the models of many physical, chemical and
biological processes has been extensively studied. In the field of
population dynamics some models of multispecies population are
described by reaction-diffusion systems. Jorne \cite{Jorne} examined
the effect of cross diffusion on the diffusive Lotka-Volterra
system. They found that the cross-diffusion may give rise to
instability in the system, although this situation seems quite rare
from an ecological point of view. Gurtin \cite{Gurtin} developed
some mathematical models for population dynamics with the inclusion
of cross-diffusion as well as self-diffusion and showed that the
effect of cross-diffusion may give rise to the segregation of two
species. Some conditions for the existence of global solutions have
been given by several authors, for example, Deuring \cite{Deuring},
Kim \cite{Kim}, Pozio and Tesei \cite{Pozio}, Yamada \cite{Yamada}.
Moreover, due to a most interesting qualitative feature: pattern
formation induced by cross-diffusion effect there are some works on
the diffusion driven instability (Turing instability \cite{Turing})
and the existence of a non-constant stationary solution, please
refer to \cite{Kuto1, Kuto2, Ling, Matano, Mimura, Peng, Wang} and
the references cited therein.

The main purpose of this paper is to study the Turing instability
which is driven solely from the effect of cross-diffusion by using
mathematical analysis and numerical simulations. The rest of this
paper is organized as follows. In section 2 we show that the unique
positive equilibrium of the ODE system (\ref{eq0}) is globally
asymptotically stable. In the section 3 we show that the positive
equilibrium remains linearly stable in the presence of
self-diffusion. It becomes linearly unstable with the inclusion of
some appropriate cross-diffusion influences. The Turing instability
occurs only when the cross-diffusion rates $k_{23}$ and $k_{32}$ are
large. The resulting patterns are computed by a numerical method and
also we devoted to some conclusions in section 4.

\section{Stability of the positive equilibrium solution of the ODE system}
In this section, we consider the stability of the positive
equilibrium solution of the system (\ref{eq0}). It is easy to know
that if
\begin{eqnarray}
\label{H} abc>\max\{e(b-a),d(a-b)\}
\end{eqnarray}
the ODE system (\ref{eq0}) has a unique positive equilibrium
$\mathbf{\bar{u}}=(\bar{u}_1,\bar{u}_2,\bar{u}_3)$ which is given by
\begin{eqnarray}
\label{U} \bar{u}_1=\frac{abc+ae-be}{abc+bd+ae},\
\bar{u}_2=\frac{abc+bd-ad}{abc+bd+ae},\
\bar{u}_3=\frac{ab(d+e)}{abc+bd+ae}.
\end{eqnarray}
We have the following result:
\begin{Theorem}\label{Th2}
The unique positive equilibrium $\mathbf{\bar{u}}$ is globally
asymptotically stable for the ODE system $(\ref{eq0})$.
\end{Theorem}
Proof. In order to prove the theorem, we need construct a Lyapunov
function for the system (\ref{eq0}).
\begin{eqnarray}
\label{Lya}
V(\mathbf{u}(t))=d(u_1-\bar{u}_1-\bar{u}_1\ln\frac{u_1}{\bar{u}_1})+e(u_2-\bar{u}_2-\bar{u}_2\ln\frac{u_2}
{\bar{u}_2})+(u_3-\bar{u}_3-\bar{u}_3\ln\frac{u_3}{\bar{u}_3}).
\end{eqnarray}
Then $V(\mathbf{\bar{u}})=0$ and $V(\mathbf{u})>0$ if
$\mathbf{u}\neq \mathbf{\bar{u}}$. By using (\ref{eq0}), we compute
\begin{eqnarray*}
\frac{\mathrm{d}V}{\mathrm{d}t}&=&d(1-\frac{\bar{u}_1}{u_1})u_1'+e(1-\frac{\bar{u}_2}{u_2})u_2'+(1-\frac{\bar{u}_3}{u_3})u_3'\\
&=&d(u_1-\bar{u}_1)[-a(u_1-\bar{u}_1)-(u_3-\bar{u}_3)]+e(u_2-\bar{u}_2)[-b(u_2-\bar{u}_2)-(u_3-\bar{u}_3)]\\
&&+(u_3-\bar{u}_3)[-c(u_3-\bar{u}_3)+d(u_1-\bar{u}_1)+e(u_2-\bar{u}_2)]\\
&=&-ad(u_1-\bar{u}_1)^2-be(u_2-\bar{u}_2)^2-c(u_3-\bar{u}_3)^2<0
\end{eqnarray*}
for all $\mathbf{u}\neq \mathbf{\bar{u}}$. By the Lyapunov-LaSalle
invariance principle \cite{Ha},   $\bar{\mathbf{u}}$ given by
(\ref{U}) is globally asymptotically stable for the kinetic
system (\ref{eq0}).

\begin{Theorem}\label{Th2}
The unique positive equilibrium $\mathbf{\bar{u}}$
 is globally asymptotically stable for the reaction-diffusion
system $(\ref{eq1})$ without cross-diffusion, i.e. $k_{ij}=0$ for
$i\neq j$.
\end{Theorem}
Proof. To study the global behavior of system (1.2), we introduce
the following Lyapunov functional
\begin{eqnarray}
\label{LP2} W(t)=\int_\Omega V(\mathbf{u}(x,t))\mathrm{d}x
\end{eqnarray}
where $V(\mathbf{u}(x,t))$ is given by (\ref{Lya}). By direct
computation, we have
\begin{eqnarray*}
\frac{\mathrm{d}W}{\mathrm{d}t}&=& \int_\Omega \mbox{grad}_{\mathbf{u}}V\cdot \frac{\partial \mathbf{u}}{\partial t}\mathrm{d}x\\
&=&\int_\Omega
\left(d(1-\frac{\bar{u}_1}{u_1}),e(1-\frac{\bar{u}_2}{u_2}),
(1-\frac{\bar{u}_3}{u_3})\right)\cdot(k_{11}\Delta u_1+u_1 f_1, k_{22}\Delta u_2+u_2 f_2,\\[1mm]
 & &k_{33}\Delta u_3+u_3 f_3)\mathrm{d}x \\ [1mm]
&=&\int_\Omega d\left(k_{11}(1-\frac{\bar{u}_1}{u_1})\Delta
u_1\right)\mathrm{d}x+\int_\Omega
e\left(k_{22}(1-\frac{\bar{u}_2}{u_2})\Delta u_2\right)\mathrm{d}x\\
& &+\int_\Omega \left(k_{33}(1-\frac{\bar{u}_3}{u_3})\Delta
u_3\right)\mathrm{d}x+\int_\Omega
\frac{\mathrm{d}V}{\mathrm{d}t}\mathrm{d}x.
\end{eqnarray*}
From Green's identity, it follows that
\begin{eqnarray*}
\int_\Omega \left(k_{ii}(1-\frac{\bar{u}_i}{u_i})\Delta u_i
\right)\mathrm{d}x&=&\int_{\partial\Omega}
k_{ii}(1-\frac{\bar{u}_i}{u_i})\frac{\partial u_i}{\partial n}
\mathrm{d}S-\int_\Omega k_{ii}\nabla _x(1-\frac{u_i}{ u_i})
\cdot\nabla_x u_i \mathrm{d}x \\
&=&-\int_\Omega k_{ii}\bar{u}_i u_i^{-2}|\nabla_x
u_i|^2\mathrm{d}x\leq 0.
\end{eqnarray*}
Since $\frac{\mathrm{d}V}{\mathrm{d}t}\leq 0$, $\int_\Omega
\frac{\mathrm{d}V}{\mathrm{d}t}\leq 0$. Thus,
$\frac{\mathrm{d}W}{\mathrm{d}t}<0$ \ for all $\mathbf{u}\neq
\bar{\mathbf{u}}$. By the Lyapunov-LaSalle invariance principle\cite{Ha},
$\bar{\mathbf{u}}$ is globally asymptotically stable for the
reaction-diffusion system (\ref{eq1}) without cross-diffusion.

\section{Effects of cross-diffusion on Turing instability}\setcounter{equation}{0}
For simplicity, we denote
\begin{displaymath} \mathbf{K}(\mathbf{u})=\left( \begin{array}{c}
(k_{11}+k_{13} u_3)u_1\\
(k_{22}+k_{23} u_3)u_2\\
(k_{31}u_1+k_{32}u_2+k_{33})u_3
\end{array} \right),\ \
  \mathbf{F}(\mathbf{u})=\left( \begin{array}{c}
au_1(1-u_1)-u_1u_3\\
bu_2(1-u_2)-u_2u_3\\
-cu_3^2+(du_1+eu_2)u_3
\end{array} \right).
\end{displaymath}
Then the reaction-diffusion system (\ref{eq1}) can be rewritten in
matrix notation as:
\begin{eqnarray}
\label{eq21}\left\{
\begin{array}{lll}
\displaystyle \frac{\partial \mathbf{u}}{\partial t}-\Delta
\mathbf{K}(\mathbf{u})=\mathbf{F}(\mathbf{u})
 \ & \mbox{in}\  \Omega\times(0,\infty),\\[2mm]
\displaystyle \frac{\partial \mathbf{u}}{\partial \eta}=0 \ & \mbox{on}\  \Omega\times(0,\infty),\\[2mm]
\mathbf{u}(x,0)=(u_{10}(x), u_{20}(x), u_{30}(x))^T \ & \mbox{in}\
\Omega.
\end{array} \right.
\end{eqnarray}
Linearizing the reaction-diffusion system (\ref{eq21}) about the
positive equilibrium $\bar{\mathbf{u}} = ( \bar{u}_1, \bar{u}_2,
\bar{u}_3)$, we have
\begin{eqnarray}
\label{L}
 \frac{\partial \Psi}{\partial t}-\mathbf{K}_{\mathbf{u}}(\mathbf{\bar{u}})\Delta \Psi=\mathbf{G}_{\mathbf{u}}(\mathbf{\bar{u}})\Psi
\end{eqnarray}
where $\Psi=(\psi_1, \psi_2, \psi_3)^T$ and
\begin{displaymath} \mathbf{K}_{\mathbf{u}}(\mathbf{u})=\left( \begin{array}{ccc}
k_{11}+k_{13} \bar{u}_3 & 0 & k_{13}\bar{u}_1\\
0 & k_{22}+k_{23} \bar{u}_3 & k_{23}\bar{u}_2\\
k_{31}\bar{u}_3 & k_{32}\bar{u}_3 &
k_{33}+k_{31}\bar{u}_1+k_{32}\bar{u}_2
\end{array} \right),
\end{displaymath}

\begin{displaymath} \mathbf{G}_{\mathbf{u}}(\mathbf{u})=\left( \begin{array}{ccc}
-a\bar{u}_1 & 0 & -\bar{u}_1\\
0 & -b\bar{u}_2 & -\bar{u}_2\\
d\bar{u}_3 & e\bar{u}_3 &-c\bar{u}_3
\end{array} \right).
\end{displaymath}
Let $0 = \mu_1<\mu_2< \mu_3< \cdot\cdot\cdot$ be the eigenvalues of
the operator $-\Delta$ on $\Omega$ with the homogeneous Neumann
boundary condition, and $E(\mu_i)$ be the eigenspace corresponding
to $\mu_i$ in $C^2(\Omega)$. Let $\mathbf{X}=\{\mathbf{u} \in
[C^1(\bar{\Omega})]^3 |\ \frac{\partial \mathbf{u}}{\partial
\eta}=0\ \mbox{on} \
\partial\Omega\}$, $\{\phi_{ij}\}_{j=1,2,...,\dim E(\mu_i)}$ be an
orthonormal basis of $E(\mu_i)$, and
$\mathbf{X}_{ij}=\{\mathbf{c}\phi_{ij}|\ \mathbf{c}\in
\mathbf{R}^3\}$. Then
\begin{eqnarray*}
\mathbf{X}=\bigoplus^\infty_{i=1}\mathbf{X}_i \ \ \mbox{and} \ \
\mathbf{X}_i=\bigoplus^{\dim E(\mu_i)}_{j=1}\mathbf{X}_{ij}.
\end{eqnarray*}
For each $i\geq 1$, $\mathbf{X}_i$ is invariant under the operator
$\mathbf{K}_{\mathbf{u}}(\mathbf{\bar{u}})\Delta +
\mathbf{G}_{\mathbf{u}}(\mathbf{\bar{u}})$. Then problem (\ref{L})
has a non-trivial solution of the form $\Psi=\mathbf{c}\phi
\exp(\lambda t)$
 if and only if $(\lambda, \mathbf{c})$ is an eigenpair for the matrix $-\mu_i \mathbf{K}_{\mathbf{u}}(\mathbf{\bar{u}})+\mathbf{G}_{\mathbf{u}}
 (\mathbf{\bar{u}})$, where $\mathbf{c}$ is a constant vector. Then the
equilibrium $\mathbf{\bar{u}}$ is unstable if at least one
eigenvalue $\lambda$ has a positive real part for some $\mu_i$.

The characteristic polynomial of $-\mu_i
\mathbf{K}_{\mathbf{u}}(\mathbf{\bar{u}})+\mathbf{G}_{\mathbf{u}}(\mathbf{\bar{u}})$
is given by
\begin{eqnarray}\label{tz}
\rho_i(\lambda)=\lambda^3+A_{2i}\lambda^2+A_{1i}\lambda+A_{0i},\label{tian1}
\end{eqnarray}
where
\begin{eqnarray}
\label{A2}
A_{2i}=(k_{11}+k_{22}+k_{33}+k_{13}\bar{u}_3+k_{23}\bar{u}_3+k_{31}\bar{u}_1+k_{32}\bar{u}_2)u_{i}+a\bar{u}_1+b\bar{u}_2+c\bar{u}_3
\end{eqnarray}
\begin{eqnarray}
\label{A1}
A_{1i}&=&[(k_{11}+k_{13}\bar{u}_3)(k_{22}+k_{23}\bar{u}_3+k_{33}+k_{13}\bar{u}_1+k_{32}\bar{u}_2)
+(k_{22}+k_{23}\bar{u}_3)\nonumber\\
&&(k_{33}+k_{31}\bar{u}_1+k_{32}\bar{u}_2)-k_{23}k_{32}\bar{u}_2\bar{u}_3-
k_{13}k_{31}\bar{u}_1\bar{u}_3]\mu_{i}^{2}\nonumber\\
&&+[(b\bar{u}_2+c\bar{u}_3)(k_{11}+k_{13}\bar{u}_3)+a\bar{u}_1(k_{22}+k_{33}+k_{23}\bar{u}_3+k_{31}\bar{u}_1+k_{32}\bar{u}_2) \ \nonumber\\
&&+b\bar{u}_2(k_{33}+k_{31}\bar{u}_1
+k_{32}\bar{u}_2)+c\bar{u}_3(k_{22}+k_{23}\bar{u}_3)+\bar{u}_2\bar{u}_3(k_{23} e-k_{32})\nonumber\\
&&+\bar{u}_1\bar{u}_3(k_{13}d-k_{31})]
\mu_{i}+a\bar{u}_1(b\bar{u}_2+c\bar{u}_3)+bc\bar{u}_2\bar{u}_3+
e\bar{u}_2\bar{u}_3+d\bar{u}_1\bar{u}_3
\end{eqnarray}
\begin{eqnarray}\label{A0}
A_{0i}&=&[(k_{33}k_{11}+k_{33}k_{13}\bar{u}_3+k_{31}k_{11}\bar{u}_1+k_{31}\bar{u}_1\bar{u}_3+k_{32}k_{11}\bar{u}_2+
k_{32}k_{13}\bar{u}_2\bar{u}_3)k_{22}\nonumber\\
&&+(k_{33}k_{11}+k_{33}k_{13}\bar{u}_3+k_{31}k_{11}\bar{u}_1)k_{23}]u_{i}^{3}+[c\bar{u}_3(k_{11}+k_{13}\bar{u}_3)
(k_{22}+k_{23}\bar{u}_3)\nonumber\\
&&+(k_{33}+k_{31}+k_{32}\bar{u}_2)[a\bar{u}_1(k_{22}+k_{23}\bar{u}_3)+b\bar{u}_2(k_{11}+k_{13}\bar{u}_3)]\nonumber\\
&&+(k_{11}+k_{13}\bar{u}_3)
(ek_{23}\bar{u}_2\bar{u}_3)-ak_{23}k_{32}\bar{u}_1\bar{u}_2\bar{u}_3+dk_{13}\bar{u}_1\bar{u}_3(k_{22}+k_{23}\bar{u}_3)\nonumber\\
&&-k_{13}k_{31}b\bar{u}_1\bar{u}_2\bar{u}_3
-k_{31}k_{22}\bar{u}_1\bar{u}_2-k_{31}k_{32}\bar{u}_1\bar{u}_3^{2}]\mu_{i}^{2}+[ac\bar{u}_1\bar{u}_3(k_{22}+k_{23}\bar{u}_3)\nonumber\\
&&+bc\bar{u}_2\bar{u}_3(k_{11}+
k_{13}\bar{u}_3)+ab\bar{u}_1\bar{u}_2(k_{33}+k_{31}\bar{u}_1+k_{32}\bar{u}_2)+e\bar{u}_2\bar{u}_3(k_{11}\nonumber\\
&&+k_{13}\bar{u}_3)+a\bar{u}_1(ek_{23}\bar{u}_2\bar{u}_3-
k_{23}\bar{u}_2\bar{u}_3)+k_{11}\bar{u}_1\bar{u}_2\bar{u}_3bd\nonumber\\
&&+d\bar{u}_1\bar{u}_3(k_{22}+k_{23}\bar{u}_3-k_{31}b\bar{u}_1\bar{u}_2\bar{u}_3)]u_{i}+(abc+ae+bd)\bar{u}_1\bar{u}_2\bar{u}_3.
\end{eqnarray}
Let $\lambda_{1i},\lambda_{2i},\lambda_{3i}$ be the three roots of
(\ref{tz}). In order to obtain the stability of $\bar{\mathbf{u}}$,
we need to show that three exists a positive constant $\delta$ such
that
\begin{eqnarray}\label{Re}
Re\{{\lambda_{1i}}\}, Re\{{\lambda_{2i}}\}, Re\{{\lambda_{3i}}\}
<-\delta, \ \mbox{for all} \ i\geq1.
\end{eqnarray}

The aim of the following theorem is to prove that the diffusion
alone (with out cross-diffusion, i.e
$k_{31}=k_{13}=k_{32}=k_{23}=0$)
   can not drive instability for this model.
 \begin{Theorem}\label{Th31} Suppose that $(\ref{H})$ holds and $k_{13}=k_{31}=k_{23}=k_{32}=0$. Then the positive equilibrium $\bar{\mathbf{u}}$ of
 $(\ref{eq21})$
 is linearly stable.
 \end{Theorem}
Proof: Substituting $k_{13}$=$k_{31}$=$k_{23}$=$k_{32}=0$  into
(\ref{A2}), (\ref{A1}) and (\ref{A0}) we have
\begin{eqnarray*}
A_{2i}&=&a\bar{u}_1+b\bar{u}_2+c\bar{u}_3+(k_{11}+k_{22}+k_{33})\mu_i>0\\
A_{1i}&=&(k_{11}k_{22}+k_{11}k_{33}+k_{22}k_{33})\mu_i^{2}+[a(k_{22}+k_{33})\bar{u}_1+b(k_{11}+
k_{33})\bar{u}_2+c(k_{11}+k_{22})\bar{u}_3]\mu_i\qquad\quad\\
&&+ab\bar{u}_1\bar{u}_2+ac\bar{u}_1\bar{u}_3+d\bar{u}_1\bar{u}_3+bc\bar{u}_2\bar{u}_3+e\bar{u}_2\bar{u}_3>0\\
A_{0i}&=&k_{11}k_{22}k_{33}\mu_i^{3}+(k_{11}k_{22}c\bar{u}_3+k_{11}k_{33}b\bar{u}_2+k_{22}k_{33}a\bar{u}_1)\mu_i^{2}\\
&&+(abk_{33}\bar{u}_1\bar{u}_2+d k_{22}\bar{u}_1\bar{u}_2+ac
k_{22}\bar{u}_1\bar{u}_3+bck_{11}\bar{u}_2\bar{u}_3+ek_{11}\bar{u}_2\bar{u}_3+bdk_{11}
\bar{u}_1\bar{u}_2\bar{u}_3)\mu_{i}\\
&&+(ae+bd+abc)\bar{u}_1\bar{u}_2\bar{u}_3>0.
\end{eqnarray*}
A direct calculation shows that $A_{2i}A_{1i}-A_{0i}>0$ for all
$i\geq1$. It follows from Routh-Hurwize criterion that all the three
roots $\lambda_{1i},\lambda_{2i},\lambda_{3i}$ of
$\rho_i(\lambda)=0$ have negative real parts for each $i\geq1$.

Let $\lambda=\mu_i\varepsilon$, then
$$\rho_i(\lambda)=\mu_i^{3}\xi^{3}+A_{2i}\mu_i^{2}\xi^{2}+A_{1i}\mu_i\xi+A_{0i}=\tilde{\rho}_i(\xi).$$
Since $\mu_i\rightarrow\infty$, as $i\rightarrow\infty$, we have
$$\bar{\rho}(\xi)=\lim_{i\rightarrow\infty}\frac{\tilde{\rho}_i(\xi)}{\mu_i^{3}}=\xi^{3}+
(k_{11}+k_{22}+k_{33})\xi^{2}+(k_{11}k_{22}+k_{22}k_{33}+
k_{11}k_{33})\xi+k_{11}k_{22}k_{33}.$$
  Applying the
Routh-Hurwitz criterion it follows that the three roots
$\xi_1,\xi_2,\xi_3$ of $\bar{\rho} (\xi)=0$ all have negative real
parts. Thus, there exists a positive constant $\bar{\delta}$ such
that $Re\{{\xi_1}\}, Re\{{\xi_2}\}$, $Re\{{\xi_3}\}\leq-
2\bar{\delta}$. By continuity, we see that there exists $i_0\geq 1$,
such that $\mu_{i0}>1$ and the three roots
$\xi_{i1},\xi_{i2},\xi_{i3}$ of $\tilde{\rho}_i(\xi) =0$ satisfy $
Re\{\xi_{i1}\},Re\{\xi_{i2}\}$, $Re\{\xi_{i3}\}
\leq-\mu_{i}\bar{\delta}\leq-\mu_{i0}\bar{\delta}\leq-\bar{\delta}$
for any $i\geq i_{0}$, Let $-\tilde{\delta}=\max_{1\leq i\leq
i_0}\{Re\{\lambda_{i1}\},Re\{\lambda_{i2}\}, Re\{\lambda_{i3}\}\} $
and $\delta=\min\{{\tilde{\delta},\bar{\delta}}\}$, Then (\ref{Re})
holds. Consequently the equilibrium $\bar{\mathbf{u}}$ is linearly
stable.

Note that $A_{2i}>0,A_{1i}>0,A_{0i}>0$, and $A_{2i}A_{1i}-A_{0i}>0$
if $k_{31}=k_{32}=0$ since the possible negative terms all involve
either $k_{31}$ or $k_{32}$. By the same arguments as in Theorem
\ref{Th31}, we have
\begin{Theorem}\label{Th3} Suppose that $(\ref{H})$ holds and $k_{31}=k_{32}=0$, Then the positive equilibrium $\bar{\mathbf{u}}$
 of $(\ref{eq1})$ is linearly stable.
\end{Theorem}

Next we consider the Turing instability i.e. the stability of the
positive equilibrium $\bar{u}=(\bar{u}_1,\bar{u}_2,\bar{u}_3)$
changing from stable, for the ODE dynamics (\ref{eq0}), to unstable
for the PDE dynamics (\ref{eq1}). Here we give sufficient conditions
on cross-diffusion which drives the instability. $k_{31}$ and
$k_{32}$ are chosen as variation parameters.
\begin{Theorem}\label{Th3}$(1)$ Suppose that $a\bar{u}_1-\bar{u}_3<0$. Consider $k_{31}$as the variation parameter, then there exists a
positive constant $ \delta_{31}$ such that when
$k_{31}>\delta_{31}$, the equilibrium $\bar{\mathbf{u}}$ is linearly unstable for some domain $\Omega$.\\
$(2)$ Suppose that $b\bar{u}_2-\bar{u}_3<0$. Consider $k_{32}$ as
the variation parameter, then there exists a positive constant
$\delta_{32}$ such that when $k_{32}>\delta_{32}$, the equilibrium
$\bar{\mathbf{u}}$ is linearly unstable for some domain $\Omega$.
\end{Theorem}
Proof: Denote
\begin{eqnarray}\label{AA}
A(\mu)=-(C_{3}\mu^{3}+C_{2}\mu^{2}+C_{1}\mu+C_{0})
\end{eqnarray}
where
\begin{eqnarray*}
C_{3}&=&[(k_{33}k_{11}+k_{33}k_{13}\bar{u}_3+k_{31}k_{11}\bar{u}_1+k_{31}\bar{u}_1\bar{u}_3+
k_{32}k_{11}\bar{u}_2+k_{32}k_{13}\bar{u}_2\bar{u}_3)k_{22}\\
&&+(k_{33}k_{11}+k_{33}k_{13}\bar{u}_3+k_{31}k_{11}\bar{u}_1)k_{23}]\\
C_{2}&=&[c\bar{u}_3(k_{11}+k_{13}\bar{u}_3)(k_{22}+k_{23}\bar{u}_3)+(k_{33}+k_{31}+k_{32}\bar{u}_2)[a\bar{u}_1(k_{22}+k_{23}\bar{u}_3)\\
&&+b\bar{u}_2(k_{11}+k_{13}\bar{u}_3)]
+(k_{11}+k_{13}\bar{u}_3)(ek_{23}\bar{u}_2\bar{u}_3)-ak_{23}k_{32}\bar{u}_1\bar{u}_2\bar{u}_3\\
&&+dk_{13}\bar{u}_1\bar{u}_3(k_{22}+k_{23}\bar{u}_3) -k_{13}k_{31}
b\bar{u}_1\bar{u}_2\bar{u}_3-k_{31}k_{22}\bar{u}_1\bar{u}_2-k_{31}k_{32}\bar{u}_1\bar{u}_3^{2}]\\
C_{1}&=&[ab\bar{u}_1\bar{u}_2k_{31}(a\bar{u}_1-\bar{u}_3)+a\bar{u}_1\bar{u}_2k_{32}(b\bar{u}_2-\bar{u}_3)+ac\bar{u}_1\bar{u}_3(k_{22}+
k_{23}\bar{u}_3)+\\
&&bc\bar{u}_2\bar{u}_3(k_{11}+k_{13}\bar{u}_3)+ab\bar{u}_1\bar{u}_2k_{33}+e\bar{u}_2\bar{u}_3(k_{11}+k_{13}\bar{u}_3)+aek_{23}\bar{u}_1
\bar{u}_2\bar{u}_3
\\&&+k_{11}bd\bar{u}_1\bar{u}_2\bar{u}_3+d\bar{u}_1\bar{u}_3(k_{22}+k_{23}\bar{u}_3)\\
C_{0}&=&(abc+ae+bd)\bar{u}_1\bar{u}_2\bar{u}_3.
\end{eqnarray*}

Case 1: $k_{31}$ is the variation parameter.\\
We assume that $a\bar{u}_1-\bar{u}_3<0$. The following arguments by
continuation are based on the fact that each root of the algebraic
equation (\ref{AA}) is a continuous function of the variation
parameter $k_{31}$. It is easy to prove that equation (\ref{AA}) has
three real roots $\mu_{1}^{(i)}=\mu_{1}^{(i)}(k_{31})$, $i=1,2,3$
when $k_{31}$ goes to infinity and they are
$\lim_{k_{31}\rightarrow\infty}\mu_{1}^{(1)}(k_{31})<0$,
$\lim_{k_{31}\rightarrow\infty}\mu_{1}^{(2)}(k_{31})=0$ and
$\lim_{k_{31}\rightarrow\infty}\mu_{1}^{(3)}(k_{31})>0$. By
continuation, there exists a positive constant $\delta_{31}$ such
that when $k_{31}>\delta_{31}$, $C_{1}>0$ and $\det(A(\mu))$ has
three real roots. Because $C_{3}>0$ and $C_{0}>0$, the mumber of
sign changes of (\ref{AA}) is exactly two. Therefore by
Descartes'rule,
the three real roots have the following properties:\\
(1) $-\infty<\mu_{1}^{(1)}<0<\mu_{1}^{(2)}<\mu_{1}^{(3)}<\infty,$\\
(2) $\det(A({\mu}))>0$ if $\mu\in(-\infty,\mu_{1}^{(1)})\cup(\mu_{1}^{(2)},\mu_{1}^{(3)}),$\\
(3) $\det(A({\mu}))<0$ if
$\mu\in(\mu_{1}^{(1)},\mu_{1}^{(2)})\cup(\mu_{1}^{(3)},\infty).$

If $\mu_{i}\in(\mu_{1}^{(2)},\mu_{1}^{(3)})$ for some $i$, then
det$(A{(\mu_{i})})>0$ by (2), and consequently
$A_{0i}=-\det(A_{(i)})<0$. The number of sign of changes of the
characteristic polynomial (\ref{tz})
$\rho_i(\lambda)=\lambda^3+A_{2i}\lambda^2+A_{1i}\lambda+A_{0i}$ is
either one or three. By Descartes'rule, the characteristic
polynomial (\ref{tz}) has at least one positive eigenvalue. Hence,
the equilibrium $\bar{\mathbf{u}}$ of (\ref{eq1}) is linearly
unstable for any domain $\Omega$ on which at least one eigenvalue
$\mu_i$ of $-\Delta$ is in the interval
$(\mu_{1}^{(2)},\mu_{1}^{(3)})$.

Case 2: $k_{32}$ is the variation parameter.\\
We assume that $b\bar{u}_2-\bar{u}_3<0$. The following arguments by
continuation are based on the fact that each root of the equation
(\ref{AA}) is a continuous function of the variation $k_{32}$. It is
easy to prove that equation (\ref{AA}) has three real roots
$\mu_{2}^{(i)}=\mu_{2}^{(i)}(k_{32})$, $i=1,2,3$ when $k_{32}$ goes
to infinity and they are
$\lim_{k_{32}\rightarrow\infty}\mu_{2}^{{(1)}}(k_{32})<0$,
$\lim_{k_{32}\rightarrow\infty}\mu_{2}^{(2)}(k_{32})=0$ and
$\lim_{k_{32}\rightarrow\infty}\mu_{2}^{(3)}(k_{32})>0$. By
continuation, there exists a positive constant $\delta_{32}$ such
that when $k_{32}>\delta_{32}$, $C_{1}>0$ and $\det(A(\mu))$
 has three real roots. Because $C_{3}>0$ and $C_{0}>0$, the mumber of sign changes of (\ref{AA}) is exactly two.
 Therefore by Descartes'rule,the three real roots have the folloeing properties:\\
(1) $-\infty<\mu_{2}^{(1)}<0<\mu_{2}^{(2)}<\mu_{2}^{(3)}<\infty$,\\
(2) $\det(A({\mu}))>0 \ \mbox{if} \ \mu\in(-\infty,\mu_{2}^{(2)})\bigcup(\mu_{2}^{(2)},\mu_{2}^{(3)})$,\\
(3) $\det(A({\mu}))<0\  \mbox{if} \
\mu\in(\mu_{2}^{(1)},\mu_{2}^{(2)})\bigcup(\mu_{2}^{(3)},\infty)$.

If $\mu_{i}\in(\mu_{2}^{(2)}, \mu_{2}^{(3)})$ for some $i$, then
$\det(A{(\mu_{i})})>0$, and consequently $A_{0i}=-\det(A_{(i)})<0$.
By similar argument as case 1, The number of sign of changes of the
characteristic polynomial (\ref{tz})
$\rho_i(\lambda)=\lambda^3+A_{2i}\lambda^2+A_{1i}\lambda+A_{0i}$ is
either one or three. By Descartes'rule, the characteristic
polynomial (\ref{tz}) has at least one positive eigenvalue. Hence,
the equilibrium $\bar{u}$ of (\ref{eq1}) is linearly unstable for
any domain $\Omega$ on which at least one eigenvalue $\mu_i$ of
$-\Delta$ is in the interval $(\mu_{2}^{(2)}, \mu_{2}^{(3)})$.

\begin{Remark}
$(i)$ In Theorem 3.3, the condition $a\bar{u}_1-\bar{u}_3<0$ and
$b\bar{u}_2-\bar{u}_3<0$ are compatible with the condition $(2.3)$
respectively.

$(ii)$ $k_{31}$ and $k_{32}$ can be chosen as variation parameters
because the number of sign of change for the polynomial $(3.8)$
could be bigger than one for large values of $k_{31}$ or $k_{32}$,
By descartes' rule, the polynomial $(3.8)$ could have positive roots
which lead to linear instability.

$(iii)$ Biological interpretation: In our model, the third species
preys on the first and second one. The positive steady state of the
model can be broken by the reaction-diffusion among two species on
the model. Case one: In this case, the first species are assumed to
reproduce exponentially unless subject to intra-species competitions
and predation. This exponential growth is represented in the
equation by the term $au_1$. The level of intra-species competitions
among the first species is assumed to be proportional to the
population density of first species by the term $au_1$. The rate of
predation upon the prey is assumed to be proportional to the rate at
which the predators and the prey meet by the term $u_1u_3$, when the
effect on first species due to the fact the third species preys on the
first one $\bar{u}_3$ are larger than the effects on first species
due to the intra-species competitions among first species
$a\bar{u}_1$, the large cross-diffusion of the third species due to
the first species $k_{31}$ can break the stability of the positive
steady state. In other words, if the predator has a dominate effect
on the decreasing of the prey such as predation rate is lager than
the rate of intra-species competitions, then the predator with large
cross-diffusion can destabilize the constant steady state. Case two:
In this case, the third species shall have a dominate effect on the
decreasing of the second species. Because $b\bar{u}_2-\bar{u}_3<0$
implies $b\bar{u}_2<\bar{u}_3$, predation rate of third species on
the second species is large than the rate of intra-species
competitions in second species. The similar situation as in case one
happens in the case two: the predator with large cross-diffusion can
destabilize the constant steady state.
\end{Remark}

\section{Numerical calculations}
In this section, using numerical methods, we illustrate that the
cross-diffusion induces spatial patterns. The initial data is taken
as
 a uniformly distributed random perturbation around the equilibrium
state $(\bar u_1, \bar u_2, \bar u_3)$ in $\Omega$, with a variance
lower than the amplitude of the final patterns. More precisely,
\bess u_{10}(x)=\bar u_1+\eta_1(x),\quad u_{20}(x)=\bar
u_2+\eta_2(x),\quad u_{30}(x)=\bar u_3+\eta_3(x), \eess where
$\eta_i\in [-1.5,1.5]$ for $i=1,2,3$. In view of Theorems \ref{Th31}
and \ref{Th3}, the Turing parameter space is (\ref{H}) under which
spatial patterns can occur.   Thus, in the system (\ref{eq1}) we fix
$a=1$, $b=1$, $c=0.1$, $d=0.1$, $e=0.1$, $k_{11}=0.1$, $k_{13}=0.1$,
$k_{22}=0.1$, $k_{23}=0.1$, $k_{31}=0.1$, $k_{33}=0.1$.
\begin{figure}[htbp]
\begin{center}
\includegraphics[ width=0.6\textwidth]{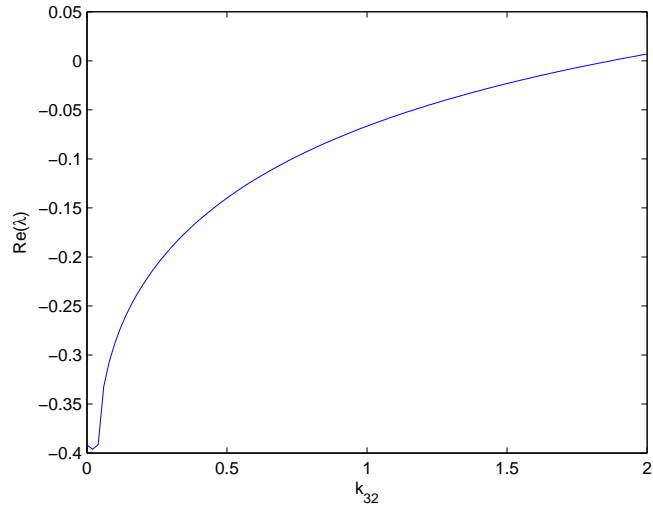}
\caption{ Dispersion relation for the real part of the eigenvalues,
$\mathrm{Re}(\lambda)$ versus the cross-diffusion coefficient
$k_{32}$. }\label{fig:fig1}
\end{center}
\end{figure}

\begin{figure}[htbp]
\begin{center}
\includegraphics[ width=0.6\textwidth]{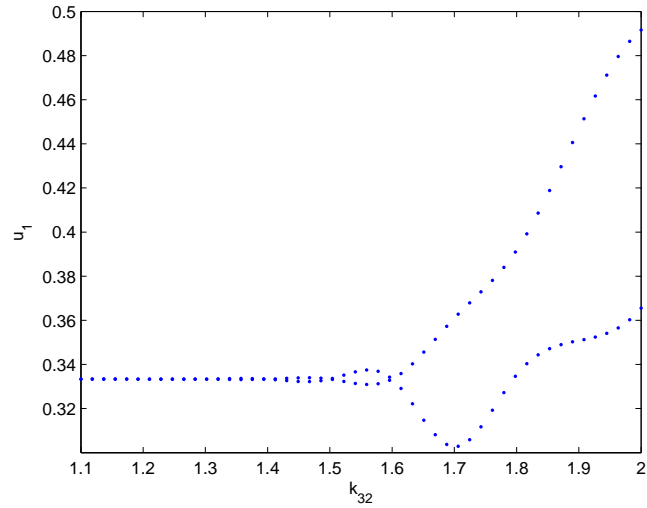}
\caption{Bifurcation diagram for Turing onset.  Maximum and minimum
$u_1$
 for different cross-diffusion in the transition from the homogeneous state to Turing pattern.  }\label{fig:fig2}
\end{center}
\end{figure}

\begin{figure}[htbp]
\begin{center}
\includegraphics[width=0.36\textwidth]{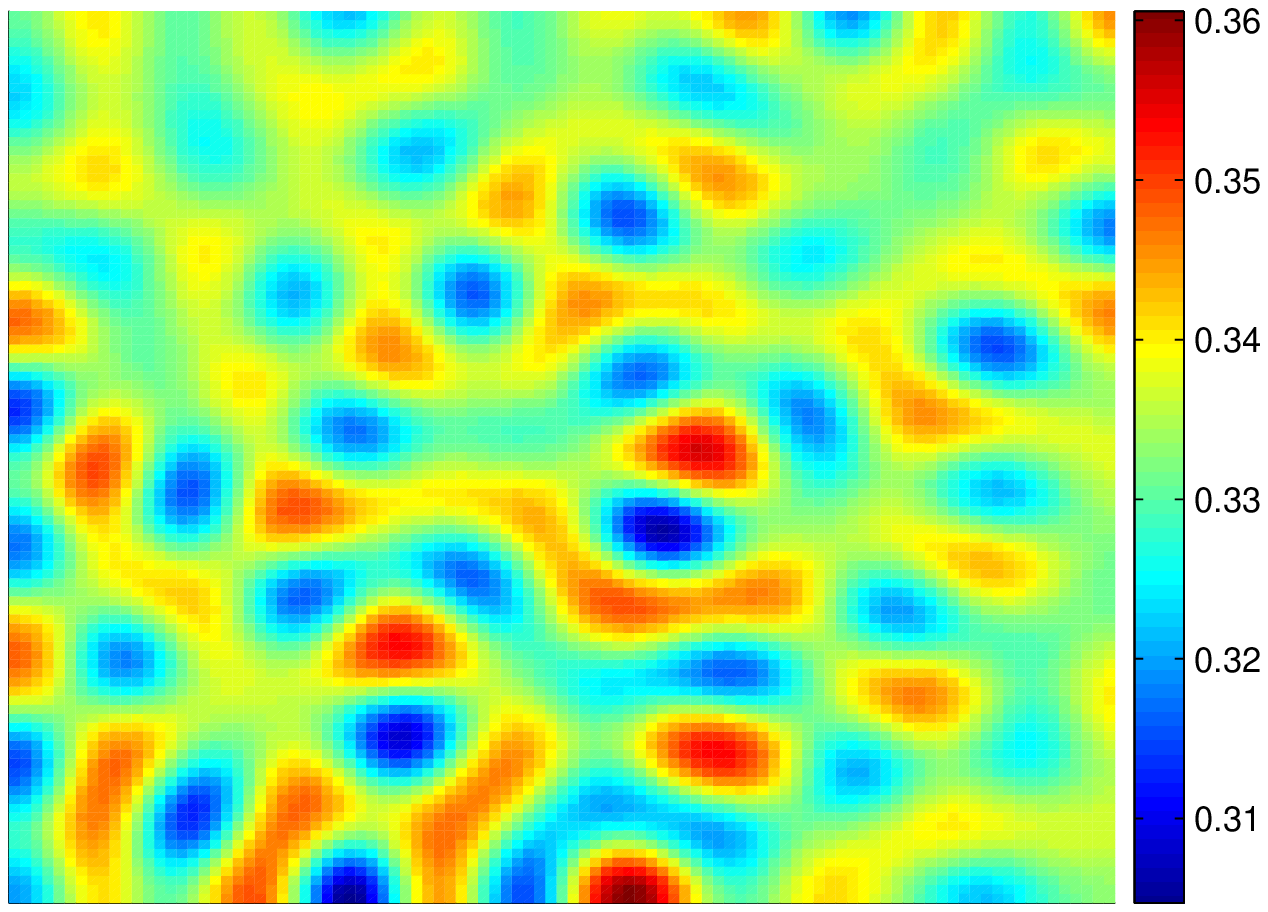}
\includegraphics[width=0.36\textwidth]{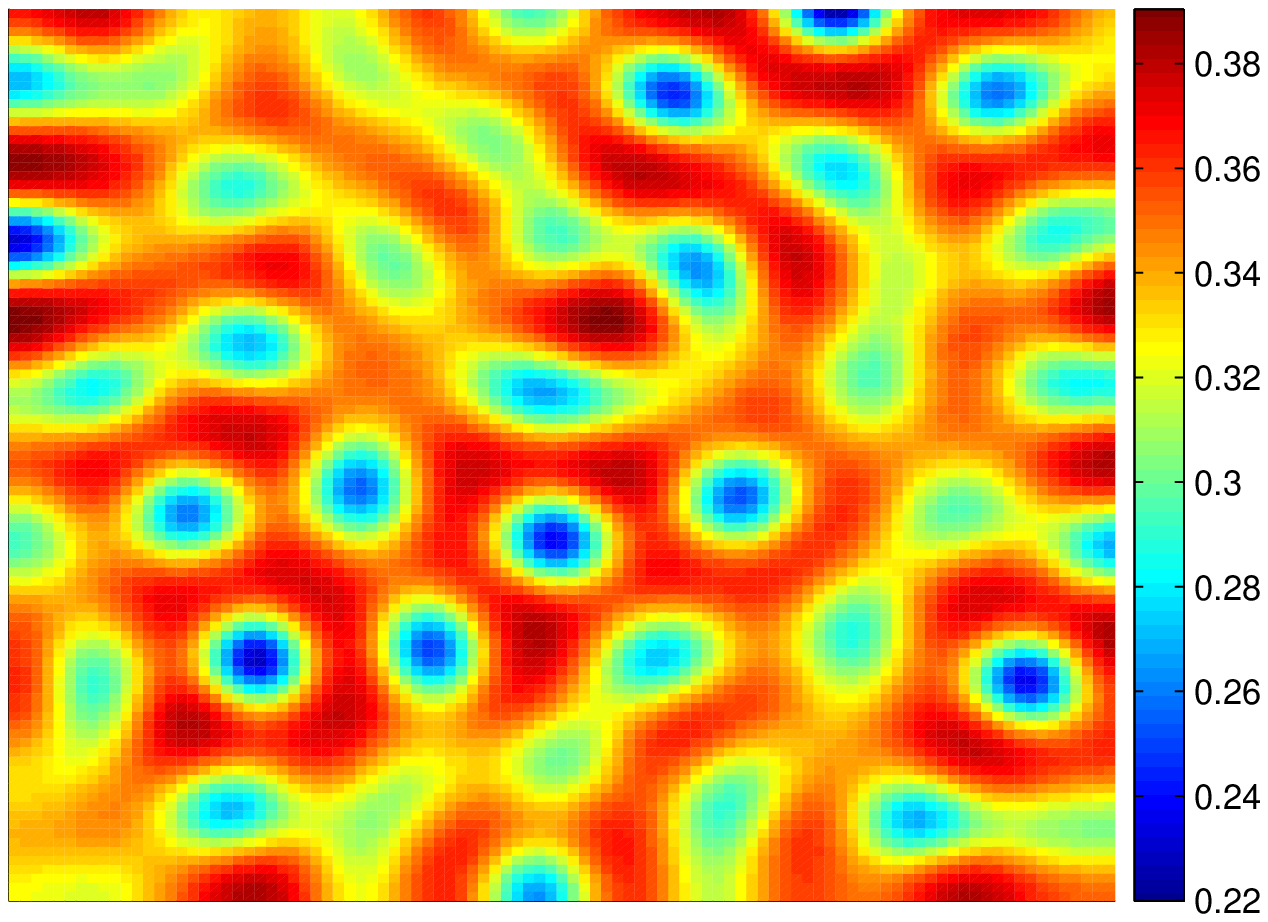}\\
\includegraphics[width=0.36\textwidth]{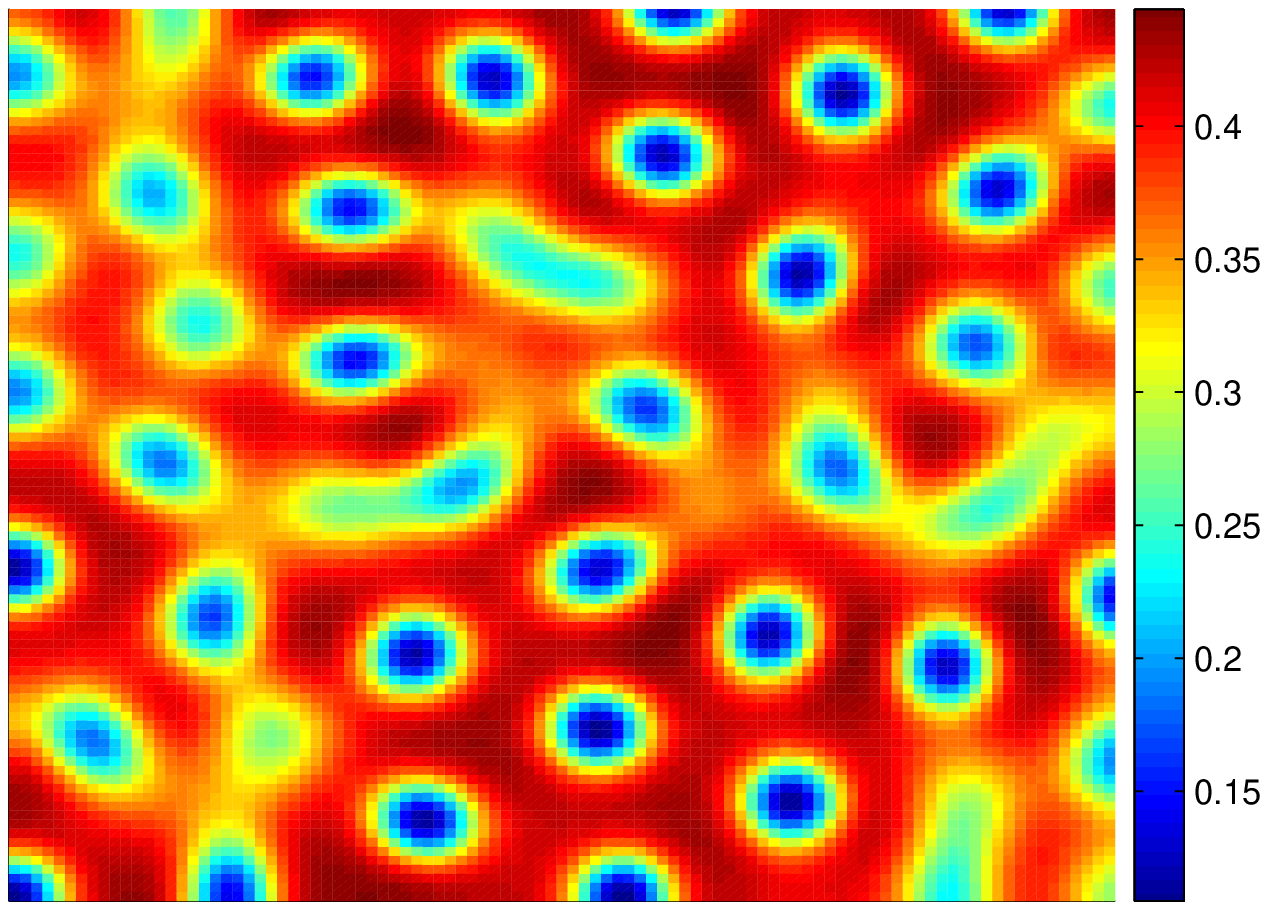}
\includegraphics[width=0.36\textwidth]{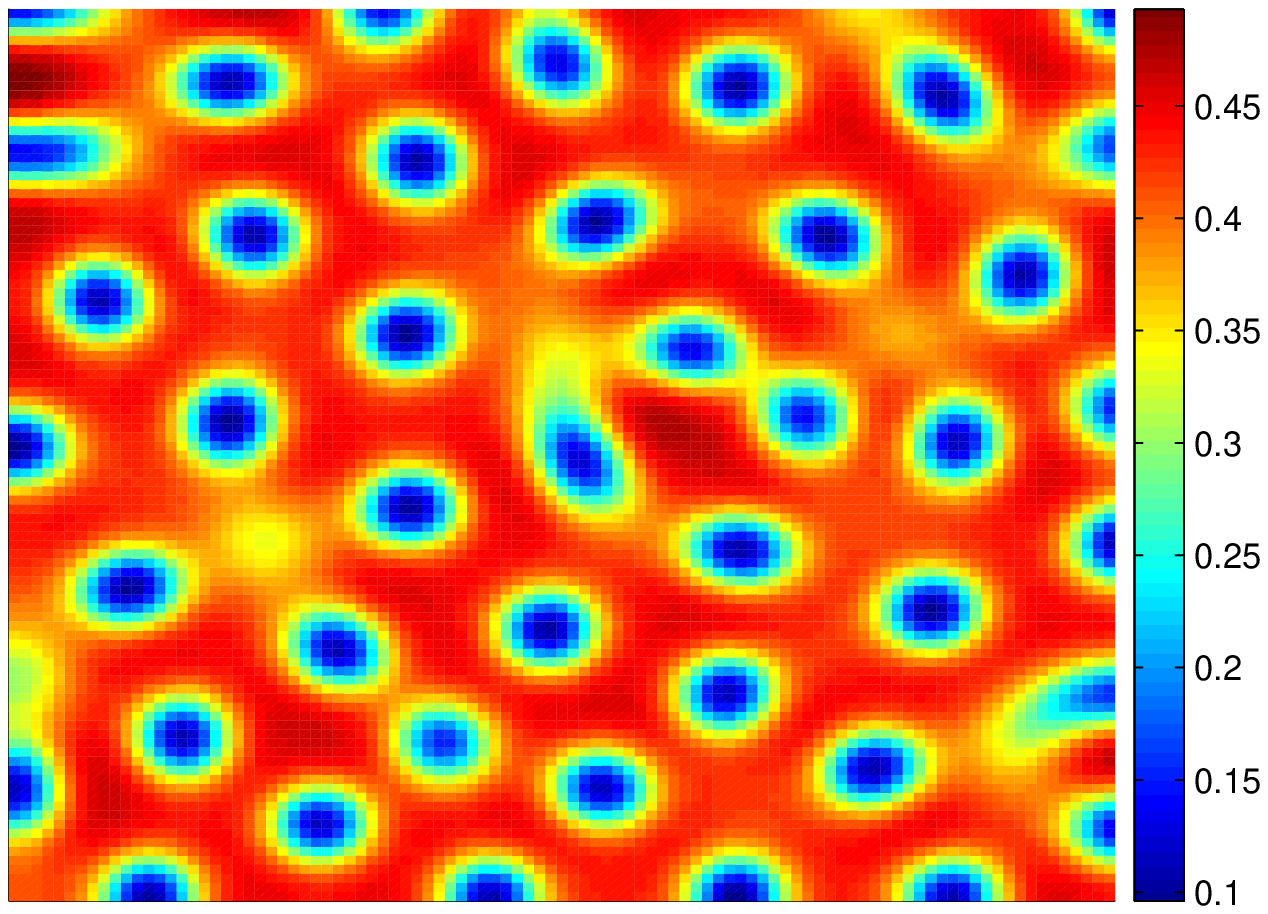}
\caption{Spatial patterns change quantitatively with different
$k_{32}$ as 1.7, 1.8, 1.9, and 2. The other parameters are stated in
the text. The steps of the iteration for time is 40000.}
\label{fig:fig3}
\end{center}
\end{figure}

In Figure \ref{fig:fig1}, we show that the real part of the
eigenvalues $\lambda$ as a function of the cross-diffusion
coefficient $d_{32}$. From the characteristic polynomial of
(\ref{tian1}), we can determine the value of $d_{32}$ such that
$Re(\lambda)>0$. Now we will implement some numerical simulations
for the system (\ref{eq1}). The domain  is confined to a square
domain $\Omega=[0, L_x]\times [0, L_y]\subset \mathbf{R}^2$. The
wavenumber for this two dimensional domain is thereby
\[\mathbf{k}=\pi(m/L_x,n/L_y), \textrm{ and }
|\mathbf{k}|=\pi\sqrt{(m/L_x)^2+(n/L_y)^2}, ~m,n=0,1,\ldots.\] We
consider system  (\ref{eq1}) in a fixed domain $L_x=40$ and
$L_y=40$, and resolve it on a grid with $100\times 100$ sites with
the space step of $\triangle x=\triangle y=1$.  For the evolution in
time, we apply a first order backward Euler time advancing scheme
with a time step $\triangle t=0.005$. By discretizing the Laplacian
in the grid with lattice sites denoted by $(i,j)$, the nine-point
formula is \bess
\triangle u|_{(i,j)}&=\frac{1}{6\triangle x^2}[4a_l(i,j)u(i-1,j)+4a_r(i,j)u(i+1,j)+4a_d(i,j)u(i,j-1)\\
&+4a_u(i,j)u(i,j+1)+a_l(i,j)u(i-1,j+1)+a_u(i,j)u(i+1,j+1)\\&+a_d(i,j)u(i-1,j-1)+a_r(i,j)u(i+1,j-1)-20u(i,j)],
\label{d1} \eess where the matrix elements of $a_l, a_r, a_d, a_u$
are unity except at the boundary. When $(i,j)$ is at the left
boundary, that is $i=0$, we define $a_l(i,j)u(i-1,j)\equiv
u(i+1,j)$, which guarantees zero-flux of reactants in the left
boundary. Similarly we define $a_r(i,j)$, $a_d(i,j)$, $a_u(i,j)$
such that the boundary is no-flux. The  nine-point formula for the
Laplacian can have a one-step error of $O(\triangle x^4)$.

In  Figure \ref{fig:fig2}, we compare the density of $u_1$ before
and after the onset of Turing patterns. Results are qualitatively
similar for $u_2$ and $u_3$, and hence omitted. In the case of
$k_{32}$ less than 1.6, i.e. the Turing instability does not occur,
we see that the density of $u_1$ is homogeneous. In the case of
$k_{32}$ larger than 1.6, i.e. Turing instability happen,  we see
that the density of $u_1$ is spatial inhomogeneous.

Now we study the change of the spatial patterns qualitatively and
quantitatively with different $k_{32}$. In general, the selection of
stripe pattern or spot pattern depends upon the non-linearities of
the reaction kinetics. Specifically, it has been shown that the
presence of quadratic nonlinearities in the reaction kinetics leads
to spot pattern, but the absence of quadratic terms leads to stripe
pattern \cite{erme}. Noticing that the reaction kinetics of
(\ref{eq1}) only has quadratic nonlinearities, in view of the theory
of pattern selection \cite{erme}, all the spatial patterns are spot
patterns. In Figure \ref{fig:fig3}, we also illustrate the
quantitative change of the spatial patterns with the different
$k_{32}$. From this simulations, we can conclude that with the
increasing of $k_{32}$, the spatial patterns converge to regular
spotted patterns. The striped patterns can not occur in our model.

\section{Comparisons and conclusions}
In this paper, we have developed a theoretical framework for
studying the phenomenon of pattern formation in a two-prey one-predator system. Applying a stability
analysis and suitable numerical simulations, we investigate the
Turing parameter space, the associated pattern type and the Turing bifurcation diagram. The proposed approach has applicability to
other reaction-diffusion systems including cross-diffusion, such as
chemotaxis and cell motility models. In this context,
it is of great interest to us the development of a general
mathematical and numerical framework that allows for the treatment
of certain degenerate quasilinear parabolic systems modeling
bacterial growth, that are known to involve several important
phenomena such as fractal morphogenesis and branching patterns.

It is worth mentioning that the authors have also the role of cross-diffusion on pattern formation for Lotka-Volterra type models in
\cite{guin,sun}. In \cite{sun}, by considering a Holling-
Tanner predator-prey model the authors investigated the Turing bifurcation and obtained the pattern selection mechanism.
In \cite{guin}, by studying the Hopf bifurcation the authors attained the spiral patterns. Apart from these work \cite{guin,sun},
what our model consider is a three species model. The difficulty is that the characteristic equation of our model
is a cubic equation. We use the continuity of the  cubic functions to overcome it. Our novelty is that we have
obtained the bifurcation diagram for Turing onset by numerical simulations, which shows the transition from the homogeneous steady state to the Turing patterns.

It is well-known that for a classical competitive model, the formation of
 patterns does not occur. We introduce the cross-diffusion into the
particular two-prey one-predator model, and show that this gives
rise to Turing-like spatial patterns. All this is confirmed with the
help of illustrating numerical simulations.







\end{document}